\documentstyle[12pt]{article}

\def\ba{\begin{eqnarray}}\def\ea{\end{eqnarray}}
\def\lb{\label}

\def\be{\begin{equation}}\def\ee{\end{equation}}

\def\R{\hat{R}}
\def\id{{I}}

\def\s{\sigma}\def\a{\alpha}
\def\d{\delta}
\def\n{\nu}\def\t{\tau}

\def\tr{{\mathrm{Tr}}\,}

\normalbaselines
\begin{document}
\begin{center}
\vspace*{1.0cm}

{\LARGE{\bf $Q$-multilinear Algebra}}

\vskip 1.5cm

{\large {\bf A. Isaev}}

\vskip 0.5 cm

Bogoliubov Laboratory of Theoretical Physics, JINR\\
141980 Dubna, Moscow Region \& Russia

\vskip 1cm

{\large {\bf O. Ogievetsky\footnote{On leave of absence from P. N.
Lebedev Physical Institute, Leninsky Pr. 53, 117924 Moscow, Russia}}}

\vskip 0.5 cm

Center of Theoretical Physics, Luminy \\
13288 Marseille \& France

\vskip 1cm

{\large {\bf P. Pyatov}}

\vskip 0.5 cm

Bogoliubov Laboratory of Theoretical Physics, JINR\\
141980 Dubna, Moscow Region \& Russia

\end{center}

\vspace{1 cm}

\begin{abstract}

The Cayley-Hamilton-Newton theorem
- which underlies the Newton identities and the Cayley-Hamilton
identity - is reviewed, first, for the classical matrices with commuting
entries, second, for two $q$-matrix algebras, the RTT-algebra and
the RLRL-algebra. The Cayley-Hamilton-Newton identities for these
$q$-algebras are related by the factorization map.
A class of algebras ${\cal{M}}(\R ,\hat{F})$
is presented. The algebras ${\cal{M}}(\R ,\hat{F})$ include the
RTT-algebra and the RLRL-algebra as particular cases. The algebra
${\cal{M}}(\R ,\hat{F})$ is defined by a pair of compatible matrices $\R$
and $\hat{F}$. The Cayley-Hamilton-Newton theorem for the algebras
${\cal{M}}(\R ,\hat{F})$ is stated. A nontrivial example of a compatible
pair is given.

\end{abstract}

\vspace{1 cm}

\section{Introduction}

In this lecture we discuss several results from the classical multilinear
algebra and their $q$-generalizations. The following simple observation
lies in the origin of our discussion: the trace of the Cayley-Hamilton
identity reproduces one of the Newton identities. Thus, it is natural
to ask if other Newton identities can be ``detraced''; in other words,
can one find matrix identities whose trace reproduces the Newton
identities. Such matrix identities do exist and we call them
the Cayley-Hamilton-Newton (CHN) identities.

In fact, we have discovered the classical CHN identities in attempts to
$q$-deform the Cayley-Hamilton theorem. We formulate the CHN identities
(which imply the Cayley-Hamilton theorem) for two standard $q$-matrix
algebras and explain the relationship of the CHN identities for
the RTT-algebra and the RLRL-algebra using the factorization map.
The formulation requires an understanding of some basic notions, like
powers of $q$-matrices and $q$-matrices acting in different copies of
a vector space. An analysis of these notions manifests different roles
of the permutation matrix and suggests a definition of a class of
$q$-matrix algebras ${\cal{M}}(\R ,\hat{F})$
which includes the standard $q$-matrix algebras. The different roles
of the permutation matrix are played by different Yang-Baxter matrices
$\R$ and $\hat{F}$ for the algebras ${\cal{M}}(\R ,\hat{F})$.

We formulate the CHN identities for the algebras ${\cal{M}}(\R ,\hat{F})$
and give a simple example.

The lecture is mainly based on papers \cite{iop1}, \cite{iops}, \cite{iop2}.

Section 2 contains an elementary introduction to the classical case.
The CHN theorem is formulated.

In Section 3 two commonly used quantum matrix algebras, the RTT-algebra
and the RLRL-algebra, are presented. A $q$-generalization of
the CHN theorem for these algebras is given.

In Section 4 we explain that in the factorizable case the CHN theorem
for the RLRL-algebra follows from the CHN theorem for the RTT-algebra.

In Section 5 we discuss more general $q$-matrix algebras for which
the $q$-CHN theorem holds. These algebras are defined by a ``compatible''
pair of Yang-Baxter matrices $\R$ and $\hat{F}$.
The section contains also a simplest nontrivial example of
a compatible pair.

\section{Classical Case}

In this lecture we shall be talking about $q$-generalizations of the
following two facts from the classical multilinear algebra:

\vskip .3cm
{\bf 2.1} Fact 1. Consider a set of $n$ variables $x_1,\dots ,x_n$.
The ring of symmetric functions in $x_1,\dots ,x_n$
has several useful sets of generators
(see, e.g. \cite{M}):

Power sums:
\be s_k:=\sum x_i^k\ ,\ k=1,\dots ,n\ ,\lb{e1}\ee
elementary symmetric functions
\be \s_k := \sum_{i_1<\dots <i_k}x_{i_1}\dots x_{i_k}\ ,
\ k=1,\dots ,n\ \lb{e2}\ee
and complete symmetric functions
\be \t_k := \sum_{i_1\leq\dots\leq i_k}x_{i_1}\dots x_{i_k}\ ,
\ k=1,\dots ,n\ .\lb{e3}\ee

The classical Newton identities relate these generating sets:
\ba
(-1)^{k+1} k\,\s_k&=&\sum_{j=1}^{k-1}(-1)^j
\ s_{k-j}\ \s_j\ ,
\lb{e4}\\[1em]
k\,\t_k&=&\sum_{j=1}^{k-1}
\ s_{k-j}\ \t_j\ .\lb{e5}
\ea
Here one additionally defines $\s_0:=1$ and $\t_0:=1$.

\vskip .3cm
{\bf 2.2} Fact 2. Let $X$ be an operator in a vector space $V$ of
dimension $n$. Let $\chi (t)$ be the characteristic polynomial of $X$,
$\chi (t) := {\mathrm det}(tI -X)$ where $I$ is the identity operator
in $V$.

The Cayley-Hamilton theorem says that $\chi (X)=0$.

\vskip .3cm
{\bf 2.3} In attempts to $q$-generalize these two facts we have found that
there is a common predecessor of both of them. To present it, we
first reformulate the Newton identities and the Cayley-Hamilton theorem.

Let $\Lambda^k V$ and ${\cal S}^k V$ be the $k$-th wedge power and
the $k$-th symmetric power of the vector space $V$ respectively.
The operator $X$ induces operators $\Lambda^k X$ and ${\cal S}^k X$ acting
in $\Lambda^k V$ and ${\cal S}^k V$: $\Lambda^k X$ and ${\cal S}^k X$
act on polyvectors as follows:
\ba
\Lambda^k X(v_1\wedge v_2\wedge \dots \wedge v_k)&:=&
Xv_1\wedge Xv_2\wedge\dots\wedge Xv_k\,\lb{e6}
\\[1em]
{\cal S}^k X(Symm(v_1\otimes v_2\otimes \dots \otimes v_k))&:=&
Symm(Xv_1\otimes Xv_2\otimes \dots \otimes Xv_k)\ .\lb{e7}
\ea

Define three sets of invariants,
\be
s_k(X):= \tr_V X^k\ ,\quad
\s_k(X):= \tr_{\Lambda^k V}(\Lambda^k X)\ , \quad
\t_k(X):= \tr_{{\cal S}^k V}({\cal S}^k X)\ .\lb{e8}
\ee

If an operator $X$ is diagonalizable with eigenvalues
$x_1,\dots ,x_n$ one finds (using a basis in which $X$
is diagonal) that $s_k(X)=s_k$, $\s_k(X)=\s_k$ and $\t_k(X)=\t_k$
where $s_k$, $\s_k$ and $\t_k$ are symmetric functions in
$x_1,\dots ,x_n$
given by (\ref{e1}), (\ref{e2}) and (\ref{e3}).

Thus, Newton identities relate traces of powers of an operator
with traces of wedge or symmetric powers of an operator.

The elementary symmetric functions $\s_k(X)$ of the operator $X$
are sums of principal $k$-minors of a matrix (in an arbitrary basis)
of the operator $X$: for a subset $\{ i_1,\dots ,i_k \} \subset
\{ 1,\dots ,\n \}$ compute the determinant of a submatrix with
rows and columns $i_1,\dots ,i_k $. The sum over all $k$-subsets
equals to $\s_k(X)$.

On the other hand the sum of principal $k$-minors is exactly what enters
the characteristic polynomial: the coefficient in $t^k$ in $\chi (t)$
gets contributions from principal $(n-k)$-minors - we have to select
$t$ on the diagonal of $tI-X$ on $k$ places, then the corresponding term
in $\chi (t) \equiv {\mathrm det}(tI -X)$  is $(-1)^{n-k}\times$
(the minor in
complementary rows and columns). In other words the characteristic
polynomial equals to
\be \sum_{j=0}^n t^j (-1)^{n-j}\s_{n-j}(X)\equiv
    (-1)^n \sum_{j=0}^n (-t)^{n-j}\s_j(X)\ .\lb{e9}\ee
Let us compare the Cayley-Hamilton theorem
\be 0=\sum_{j=0}^n (-X)^{n-j}\s_j(X)\ \lb{e10}\ee
with the Newton identities (\ref{e4})
\be k\s_k(X)=-\sum_{j=0}^{k-1}(-1)^{k-j}s_{k-j}(X)\s_j(X)\equiv
   -\sum_{j=0}^{k-1}{\mathrm Tr}(-X)^{k-j}\s_j(X)\ .\lb{e11}\ee
The comparison makes clear that the trace of the Cayley-Hamilton
identity (\ref{e10}) reproduces the Newton identity with $k=n+1$
($\s_{n+1}=0$, the $(n+1)$-st wedge power of an $n$-dimensional
space vanishes).

Therefore we conclude that there is a natural way to ``detrace''
the $(n+1)$-st Newton identity. A question arises: can we detrace
other Newton identities?

\vskip .3cm
{\bf 2.4} It turns out that the answer is positive.

\vskip .3cm
{\bf Theorem.} Let $X^{[j]}:={\mathrm Tr}_{(1\dots j-1)}
\Lambda^jX$ (the operator $\Lambda^jX$ acts on $j$-polyvectors;
we take the trace in all indices but the last one). Then
\be kX^{[k]}=-\sum_{j=0}^{k-1}(-X)^{k-j}\s_j(X) \ \ \ \forall \ k \ .
\lb{e12}\ee

\vskip .3cm
This theorem underlies both Newton and Cayley-Hamilton identities:
the trace of eqs. (\ref{e12}) gives the Newton identities;
eq. (\ref{e12}) for $k=n+1$ is the Cayley-Hamilton identity.
We couldn't find this theorem in the literature and we gave
it a name: the Cayley-Hamilton-Newton (CHN) theorem.

\vskip .3cm
{\bf 2.5} {\it Remarks.} 1. There is the Cayley-Hamilton-Newton theorem
for the symmetric powers. Let $X^{(j)}:={\mathrm Tr}_{(1\dots j-1)}
S^jX$. Then
\be kX^{(k)}=\sum_{j=0}^{k-1}X^{k-j}\t_j(X) \ \ \ \forall \ k \ .
\lb{e13}\ee
Now the sequence $X^{(j)}$ does not terminate; there is no analogue
of the Cayley-Hamilton theorem. However, taking the trace of
(\ref{e13}) one finds another sequence (\ref{e5}) of the Newton
identities.

\vskip .3cm
2. A matrix $X$ generates an algebra $k[X]$ of polynomials
in $X$. As a vector space, $k[X]$ has a filtration
$F_jk[X]=<1,X,\dots ,X^j>$ which stabilizes after at most $n$ steps.
The CHN theorem implies that $X^{[j]}$ (and $X^{(j)}$) belongs to
$k[X]$ and moreover to $F_jk[X]$.

\vskip .3cm
3. The following qualitative argument shows that an identity
of the type (\ref{e12}) should exist, or, in other words,
that the operator $X^{[k]}$ should be expressible in terms of
usual powers $X^j$.

The space $\Lambda^k V$ can be embedded into $V^{\otimes k}$ as
an image of a projector (antisymmetrizer) $A_k$ acting in
$V^{\otimes k}$,
\be (A_k)^{i_1\dots i_k}_{j_1\dots j_k}=\frac{1}{k!}
    \delta^{i_1}_{[j_1}\dots \delta^{i_k}_{j_k]} \ .\lb{e14}\ee
Here $[\dots ]$ means antisymmetrization.

The space $V^{\otimes k}$ decomposes into a direct sum of
subspaces corresponding to the Young diagrams. Therefore, the
subspace $\Lambda^k V$ has a well defined complement. One can
interpret the operator $\Lambda^k X$ as an operator acting in
$V^{\otimes k}$ as $\Lambda^k X = A_k X\otimes \dots\otimes X$
($k$ times): its restriction to $\Lambda^k V$ coincides with
(\ref{e6}) and it acts as zero on the complement to $\Lambda^k V$.

Eq. (\ref{e14}) shows that $\Lambda^k X$ equals
\be (\Lambda^kX)^{i_1\dots i_k}_{j_1\dots j_k}=\frac{1}{k!}
    X^{i_1}_{[j_1}\dots X^{i_k}_{j_k]} \ .\lb{e15}\ee
This is a sum of $k!$ terms. Tracing each term in the spaces
$1,\dots ,k-1$ one obtains some power of $X$ times some products of
traces of powers of $X$. Therefore we have
$X^{[k]}=\sum X^{k-j}\mu_j(X)$. Eq. (12) gives an exact expression for
the scalar coefficients $\mu_j(X)$.

Similar arguments hold for the CHN identities (\ref{e13}).

\section{Two quantum matrix algebras}

{\bf 3.1} The standard Drinfeld-Jimbo $\hat{R}$-matrix
\be \R^{ij}_{kl}=q^{\d_{ij}}\ \d^i_l\d^j_k +(q-q^{-1})\theta (l-k)
    \d^i_k\d^j_l \lb{e16}\ee
has two eigenvalues, $q$ and $(-q^{-1})$. In other words, its
projector decomposition has two terms, $\R =qS-q^{-1}A$. Conventionally,
the projector $S$ is called $q$-symmetrizer and the projector $A$ -
$q$-antisymmetrizer.

One says that Yang-Baxter matrices having two projectors are of Hecke
type. The $q$-CHN theorem which we shall formulate is valid for
arbitrary Yang-Baxter matrices, without any conditions on ranks of
projectors $S$ and $A$.

\vskip .3cm
{\bf 3.2} Two types of $q$-matrix algebras are commonly used in
the literature.

The first one is the RTT-algebra \cite{FRT}. It is generated by a matrix
$T^i_j$ and relations
\be \R T_1T_2=T_1T_2\R\ .\lb{e17}\ee
Here $T_1$ is the matrix $T$ in the first copy of the space $V$,
$T_2$ is the matrix $T$ in the second copy of the space $V$.

This algebra is the algebra of functions on the $q$-group.

\vskip .3cm
The second algebra is the RLRL-algebra (see \cite{KS} and
references therein).
It is generated by a matrix $L^i_j$ and relations
\be \R L_1\R L_1=L_1\R L_1\R\ .\lb{e18}\ee
For a deformation-type Hecke $\R$ (like $\R$ in (\ref{e16})) this algebra
can be interpreted as the $q$-universal enveloping algebra: in the
first orders, $\R =P(1+\a r)+O(\a^2)$ and $L=1+\a l+O(\a^2)$ ($P$ is the
permutation matrix and $\a$ is the deformation parameter); the eq. (\ref{e18})
implies, in the second order in $\a$, that $[l_1,l_2]=[r+r_{21},l_1]$;
the sum $r+r_{21}$ is proportional to the permutation $P$; the commutation
relations $[l_1,l_2]=[P,l_1]$ is just a compact way to write the
commutation relations for the Lie algebra $gl(n)$ (or, with
a condition of zero traces, for $sl(n)$).

\vskip .3cm
{\bf 3.3} {\bf CHN theorem for the RTT-algebra}. The result
will look almost the
same as in the classical case, we only have to explain meanings
of $T^{[k]}$, $T^k$, $\s_k(T)$ (in the present talk we shall give
a $q$-version of the
eq. (\ref{e12}) only; for a $q$-version of (\ref{e13})
see \cite{iop1}).

In the Hecke situation there is a well defined sequence of projectors
(we assume that $q$ is not a root of unity), antisymmetrizers,
defined inductively \cite{G} by
\be A_1:=\id\  ,\ A_k:={1\over k_q}\,
A_{k{-}1}\left(q^{k-1}-(k{-}1)_q\R_{k{-}1}\right)A_{k{-}1}\ .
\lb{e19}\ee
Here $\R_{k{-}1}$ is the operator $\R$ acting in the $(k-1)$-st and
$k$-th copies of the space $V$ and $k_q$ is the $q$-number,
$k_q:=(q^k-q^{-k})/(q-q^{-1})$. The second antisymmetrizer, $A_2$,
coincides with the projector $A$ entering the spectral decomposition
of $\R$.

To define $X^{[k]}$ classically, it does not matter if we take traces
in the spaces $1,\dots ,k-1$ or $2,\dots ,k$. For a general $\R$
these two possibilities differ and we define
two versions of the $k$-th wedge power of the matrix $T$:
\be T^{\underline{[k]}} := \tr_{(1 \dots k-1)}
\left(A_{k} T_1  \dots T_k  \right) \ \lb{e20}\ee
and
\be T^{\overline{[k]}} := \tr_{(2 \dots k)}
\left( A_{k} T_1  \dots T_k \right)\ .\lb{e21}\ee
The elementary symmetric functions in the ``spectrum of $T$''
are defined by
\be \s_k(T):=q^k \tr_{(1 \dots k)}
\left(A_{k} T_1  \dots T_k  \right) \ .\lb{e22}\ee

The definition of powers of the matrix $T$ is more interesting.
To $q$-deform a classical object or notion whose
definition involves the permutation matrix, one always has to
analyse whether the permutation stays a permutation on the
$q$-level or it becomes the $\R$-matrix. There is a way to define
the square of a classical matrix $X$ using the permutation
matrix $P$: $X^2={\mathrm tr}_1(PX_1X_2)$. It turns out that the
right choice is to replace $P$ by $\R$. Again, there are
two versions:
\be T^{\underline{k}} := \tr_{(1 \dots k-1)}
\left(\R_1 \R_2 \dots \R_{k-1} T_1 T_2 \dots T_k  \right) \
\lb{e23}\ee
and
\be T^{\overline{k}} := \tr_{(2 \dots k)}
\left( \R_1 \R_2 \dots \R_{k-1} T_1 T_2 \dots T_k \right)\ .
\lb{e24}\ee
With these preliminaries we can formulate the $q$-CHN theorem
for the RTT-algebra:

\vskip .3cm
{\bf Theorem.} The following identities hold:
\be k_q\, T^{\underline{[k]}} =
  -\sum_{j=0}^{k-1}(-1)^{k-j}\s_j(T)\, T^{\underline{k-j}}
  \ \ \ \forall\ k\
\lb{e25}\ee
and
\be k_q\, T^{\overline{[k]}} =
  -\sum_{j=0}^{k-1}(-1)^{k-j} T^{\overline{k-j}}\, \s_j(T)
  \ \ \ \forall\ k\ .\lb{e26}\ee

\vskip .3cm
{\it Remarks.} 1. The elementary symmetric functions $\s_j(T)$
form a commutative set but they are not central; the order
of terms in the right hand sides of (\ref{e25}) and (\ref{e26})
is essential.

\vskip .3cm
2. Taking $k=n+1$ one obtains two $q$-versions of the characteristic
identity. However, there is only one version of the Newton identities:
taking trace of either (\ref{e25}) or (\ref{e26}) one finds that
$s_k(T)$ can be expressed in terms of $\s_j(T)$'s. Therefore
$s_k(T)$ commute with $\s_j(T)$ and the tracing of the identities
(\ref{e25}) and (\ref{e26}) produce the same result.

\vskip .3cm
{\bf 3.4} {\bf CHN theorem for the RLRL-algebra}. Again, the
result looks as the
classical one after we give the right meaning to the notation.

It turns out that for the RLRL-algebra the subtle point is the
definition of the $L$-matrix ``acting in a $k$-th copy'' of the space
$V$. It is not at all $L_k$ and moreover, it acts in all spaces
from 1 till $k$. We denote it by $L_{\overline{k}}$ (to distinguish
from $L_k$). The definition is inductive:
\be  L_{\overline{1}}:=L_1\ ,\
     L_{\overline{k}}:=\R_{k-1}L_{\overline{k}}\R_{k-1}^{-1}\ .
\lb{e27}\ee

The definitions of $L^{[k]}$ and $\s_j(L)$ are as follows:
\be L^{[k]}:={{\mathrm Tr}_q}_{(2\dots k)} (A_{k}L_{\overline{1}}\dots
    L_{\overline{k}})\ \lb{e28}\ee
and
\be \s_j(L):={{\mathrm Tr}_q}_{(1\dots j)} (A_{j}L_{\overline{1}}\dots
    L_{\overline{j}})\ .\lb{e29}\ee
Here ${\mathrm Tr}_q$ is the $q$-trace, ${\mathrm Tr}_q(Z):=
{\mathrm Tr}(DZ)$ for an arbitrary matrix $Z$; the matrix $D$
is defined by
\be {\mathrm Tr}_2(\R D)=I\ .\lb{e30}\ee

Finally, the $q$-power of the matrix $L$ is just the usual
power.

We are ready to state the $q$-CHN theorem for the RLRL-algebra.

\vskip .3cm
{\bf Theorem.} The following identities hold:
\be k_q\, L^{[k]} =
  -\sum_{j=0}^{k-1}\s_j(L)\ (-L)^{k-j}
  \ \ \ \forall\ k\ .\lb{e31}\ee

\vskip .3cm
{\it Remark.} The elements $\s_j(L)$ are central; there is only
one version of the CHN identities.

\section{Factorization}

Denote by ${\cal{U}}$ the algebra dual to the RTT-algebra and
by ${\cal{U}}^*$ the RTT-algebra itself. Assume that ${\cal{U}}$
is quasitriangular with the universal $R$-matrix
${\cal{R}}\in {\cal{U}}\otimes {\cal{U}}$. Assume also that
the numerical matrix $P\R$ is the image of ${\cal{R}}$ in the
representation in the vector space $V$.

Contracting the second argument of the element
${\cal{R}}_{21}{\cal{R}}\in {\cal{U}}\otimes {\cal{U}}$ with
an arbitrary element $x\in {\cal{U}}^*$ we obtain an element
from ${\cal{U}}$. This defines a mapping $\phi :
{\cal{U}}^*\rightarrow {\cal{U}}$ which is called the
factorization map (\cite{RS},\cite{D}),
\be \phi (x):=<{\cal{R}}_{21}{\cal{R}},x>_2\ \ \ ,\
    x\in {\cal{U}}^*\ .\lb{e32}\ee
Here $<,>$ is the pairing between ${\cal{U}}$ and ${\cal{U}}^*$,
the index 2 in $<,>_2$ means the second copy of ${\cal{U}}$ of
the first argument.

The map $\phi$ is not a homomorphism, the matrix $L=\phi (T)$
satisfies the RLRL-algebra (\ref{e18}). However, as a linear
map, $\phi$ transforms any identity to an identity.

The matrix $D$ defined by (\ref{e30}) satisfies the
equality $\R D_1D_2=D_1D_2\R$. It follows then that
the matrix $\tilde{T}:=DT$ satisfies (\ref{e17}) or, in
other words, $T\mapsto \tilde{T}$ is an automorphism of
the RTT-algebra. Since $\tilde{T}$ generates the same
RTT-algebra, the matrix $\tilde{T}$ satisfies the CHN
identities (\ref{e26}).

\vskip .3cm
{\bf Theorem.} The map $\phi$ transforms the CHN identities
for $\tilde{T}$ (\ref{e26}) to the CHN identities (\ref{e31})
for $L=\phi (T)$.

\vskip .3cm
The algebra ${\cal{U}}$ is called factorizable if the map
$\phi$ is an isomorphism of the underlying vector spaces.

For a factorizable ${\cal{U}}$ the theorem above gives
another way to prove the CHN identities (\ref{e31}) for
the RLRL-algebra using the CHN identities (\ref{e26}) for
the RTT-algebra.

\section{Further $q$-generalizations}

{\bf 5.1} The statement of the CHN theorem for the RTT-algebra and the
RLRL-algebra looks the same but the meaning of objects
(like powers of matrices, wedge powers of matrices,
matrices acting in a $j$-th copy of the space) is different.
It seems natural to try to construct a wider class of algebras
among which the RTT-algebra and the RLRL-algebra are just particular
cases.

Such a class of algebras indeed exists. The idea behind
the construction appears already at the classical level
and again shows up in the analysis of the roles of the
permutation matrix. There are two quite different uses
of the permutation matrix $P$. Let $X$ be a matrix acting in the
vector space $V$. Consider the tensor square $V\otimes V$
of the space $V$ and let $X_1$ be the operator $X$ acting
in the first copy of $V$, $X_1=X\otimes I$. The operator
$X_2$ acting in the second copy of $V$ can be obtained
from $X_1$ with the help of the permutation matrix:
$X_2=PX_1P$. This is the first use of $P$.

The second use of $P$: to say that the matrix elements of $X$
commute, one writes $PX_1X_2=X_1X_2P$.

At the $q$-level we can play with different possibilities
of $q$-deforming the uses of $P$. In particular, there is a room
for two different Yang-Baxter matrices to appear and one can
expect an existence of $q$-algebras defined by a pair of
Yang-Baxter matrices $\R$ and $\hat{F}$, the matrix $\R$
governs the commutation relations between the matrix elements
of a quantum matrix, the matrix $\hat{F}$ is responsible
for shifting from a copy of the space $V$ to the next copy.

A detailed analysis shows that the Yang-Baxter matrices $\R$ and
$\hat{F}$ should be ``compatible'' in the following
sense:
\be \R_1\hat{F}_2\hat{F}_1=\hat{F}_2\hat{F}_1\R_2\ ,\
    \R_2\hat{F}_1\hat{F}_2=\hat{F}_1\hat{F}_2\R_1\ .\lb{e33}\ee
Such compatible pairs appear in a particular kind \cite{R} of
the Drinfeld twist \cite{D2}, the matrix $\R^F=
\hat{F}\R\hat{F}^{-1}$ is again a Yang-Baxter matrix.

One easily checks that the pair $\R^F$ and $\hat{F}$ is again
compatible and one can twist the second time to obtain a
Yang-Baxter matrix $\R^{FF}$.

The generalized algebra ${\cal{M}}(\R ,\hat{F})$ is the
algebra generated by a matrix $M^i_j$ and relations
\be \R M_{\overline 1}M_{\overline 2}=
    M_{\overline 1}M_{\overline 2}\R^{FF}\ .\lb{e34}\ee
Here $M_{\overline 1}:=M_1$ and $M_{\overline 2}:=
\hat{F}M_{\overline 1}\hat{F}^{-1}$.

\vskip .1cm
For $\hat{F}=P$ the algebra ${\cal{M}}(\R ,\hat{F})$ becomes
the RTT-algebra; for $\hat{F}=\R$ the algebra
${\cal{M}}(\R ,\hat{F})$ becomes the RLRL-algebra.

{\bf 5.2} To conclude, we shall formulate the $q$-CHN theorem
for the algebra ${\cal{M}}(\R ,\hat{F})$. To this end,
define the $k$-th power of $M$ by
\be M^{\overline{k}} := \tr_{F(2 \dots k)}
\left( \R_1 \R_2 \dots \R_{k-1}
M_{\overline{1}} M_{\overline{2}} \dots
M_{\overline{k}} \right)\ ,
\lb{e35}\ee
the $k$-th wedge power of $M$ by
\be M^{\overline{[k]}} := \tr_{F(2 \dots k)}
\left( A_{k} M_{\overline{1}}  \dots
M_{\overline{k}} \right)\ \lb{e36}\ee
and the elementary symmetric functions in the spectrum of $M$
by
\be \s_k(M):={\mathrm Tr}_F (M^{[k]})\ .\lb{e37}\ee
Here the matrices $M_{\overline{k}}$ are defined
inductively, $M_{\overline{k+1}}:=\hat{F}_k
M_{\overline{k}}\hat{F}_k^{-1}$; the
antisymmetrizer $A_{k}$ is built with the help of the
Yang-Baxter matrix $\R$ by eqs. (\ref{e19}); $\tr_{F}$ is the quantum
trace defined by the Yang-Baxter matrix $F$, that is,
for an arbitrary matrix $Z$, ${\mathrm Tr}_F(Z):=
{\mathrm Tr}(D(F)Z)$ where the matrix $D(F)$
is defined by ${\mathrm Tr}_2(\hat{F} D(F))=I$.

\vskip .3cm
{\bf Theorem}. The following identities hold
\be (-1)^{k-1} k_q\ M^{\overline{[k]}} =
  \sum_{j=0}^{k-1}(-q)^{j} M^{\overline{k-j}}\, \s_j(M)
  \ \ \ \forall\ k\ .\lb{e38}\ee

\vskip .3cm
There are versions of the $q$-CHN theorem for
$M^{\underline{[k]}}$ and for the $q$-symmetric powers of $M$
as well \cite{iop2}.

\vskip .3cm
{\bf 5.3} For the standard Yang-Baxter matrix (\ref{e16}) the
multiparametric deformation is given by a twist with the
Yang-Baxter matrix $\hat{F}=\Delta P$, where $\Delta$ is a
diagonal matrix. In this case one has $\R^{FF}=\R$.

We shall give a simplest example of a compatible pair of Yang-Baxter matrices
$\R$ and $\hat{F}$, for which $\R^{FF}\neq\R$.

Consider the Cremmer-Gervais $R$-matrix \cite{CG},
\be \R (b,y)  = \! { \left( \begin{array}{ccccccccc}
q & \cdot & \cdot \, & \, \cdot & \cdot & \cdot \, & \, \cdot &
\cdot & \cdot \\
\cdot & \cdot & \cdot \, & \, b & \cdot & \cdot \, & \, \cdot &
\cdot & \cdot \\
\cdot & \cdot & \cdot \, & \, \cdot & \cdot & \cdot \, & \, b^2/q &
\cdot & \cdot\\[1mm]
\cdot & 1/b & \cdot \, & \, \lambda & \cdot & \cdot \, & \, \cdot &
\cdot & \cdot \\
\cdot & \cdot & y \, & \, \cdot & q & \cdot \, & \, - b^2 y /q^2 &
\cdot & \cdot \\
\cdot & \cdot & \cdot \, & \, \cdot & \cdot & \cdot \, & \, \cdot &
b & \cdot\\[1mm]
\cdot & \cdot & q/b^2 \, & \, \cdot & \cdot & \cdot \, & \, \lambda &
\cdot & \cdot \\
\cdot & \cdot & \cdot \, & \, \cdot & \cdot & 1/b \, & \, \cdot &
\lambda & \cdot \\
\cdot & \cdot & \cdot \, & \, \cdot & \cdot & \cdot \, & \, \cdot &
\cdot & q
\end{array} \right) }\ \ .
\lb{e39}\ee
Here $\lambda=q-1/q$; zero entries of the matrix are denoted by dots.

The parameter $y$ is irrelevant, it can be set to 1 by a rescaling of the
coordinates of the space $V$. The parameter $y$ is introduced for
convenience, with its help we will write the matrix $\R^{FF}$ in a
compact way.

The matrix $\hat{F}$ has the form
\be \hat{F}  = \! { \left( \begin{array}{ccccccccc}
1 & \cdot & \cdot \, & \, \cdot & \cdot & \cdot \, & \, \cdot &
\cdot & \cdot \\
\cdot & \cdot & \cdot \, & \, \beta & \cdot & \cdot \, & \, \cdot &
\cdot & \cdot \\
\cdot & \cdot & \cdot \, & \, \cdot & \cdot & \cdot \, & \, -1 &
\cdot & \cdot\\[1mm]
\cdot & \alpha & \cdot \, & \, \cdot & \cdot & \cdot \, & \, \cdot &
\cdot & \cdot \\
\cdot & \cdot & \cdot \, & \, \cdot & \gamma & \cdot \, & \, \cdot &
\cdot & \cdot \\
\cdot & \cdot & \cdot \, & \, \cdot & \cdot & \cdot \, & \, \cdot &
\alpha & \cdot\\[1mm]
\cdot & \cdot & -1 \, & \, \cdot & \cdot & \cdot \, & \, \cdot &
\cdot & \cdot \\
\cdot & \cdot & \cdot \, & \, \cdot & \cdot & \beta \, & \, \cdot &
\cdot & \cdot \\
\cdot & \cdot & \cdot \, & \, \cdot & \cdot & \cdot \, & \, \cdot &
\cdot & 1
\end{array} \right) } \ \ , \;\lb{e40}\ee
where $\alpha^2=\beta^2=\gamma^2=-1$.

The Yang-Baxter matrices (\ref{e39}) and (\ref{e40}) form a compatible pair.

In \cite{iop2} it is shown that for a compatible pair $\R$
and $\hat{F}$ one has
\be \R^{FF}=D(F)_1D(F)_2\ \R \ (D(F)_1D(F)_2)^{-1}\ .\lb{e41}
\ee

The matrix $D(F)$ for the Yang-Baxter matrix $\hat{F}$ is
given by
\be D(F)=diag\{1,\gamma^{-1},1\}\ .\lb{e42}\ee

Finally, for the compatible pair (\ref{e39}) and (\ref{e40})
one finds
\be D(F)_1 D(F)_2\ R(b,y)\ (D(F)_1 D(F)_2)^{-1} = R(b,-y)\ .
\lb{e43}\ee

\section*{Acknowledgments}

This work was supported in parts by CNRS grant PICS 608
and RFBR grant 98-01-22033. The work of AI and PP was also supported by
RFBR grant 97-01-01041.

\end{document}